\input amstex
\documentstyle{amsppt}
\magnification 1200
\NoRunningHeads
\NoBlackBoxes
\document

\def\y{\bold y}
\def\Res{\text{Res}}
\def\v{{\bold v}}
\def\w{{\bold w}}

\def\ell{{\text{ell}}}

\def\u{\bold u}

\def\qdet{\text{qdet}}

\def\RR{\Bbb R}

\def\h1{\hat{\bold 1}}

\def\Hom{\text{Hom}}

\def\Ua{U_q(\tilde\g)}
\def\U2{{\Ua}_2}
\def\g{\frak g}

\def\Z{\Bbb Z}
\def\C{\Bbb C}

\def\z{\bold z}

\def\<{\langle}
\def\>{\rangle}
\def\o{\otimes}
\def\e{\varepsilon}

\def\End{\text{End}}

\def\h{{\frak h}}
\topmatter
\title Quantization of Lie bialgebras, V
\endtitle
\author {\rm {\bf Pavel Etingof and David Kazhdan} \linebreak
\vskip .1in
Department of Mathematics\linebreak
Harvard University\linebreak 
Cambridge, MA 02138, USA\linebreak
e-mail: etingof\@math.harvard.edu\linebreak kazhdan\@math.harvard.edu}
\endauthor
\endtopmatter

This paper is a continuation of \cite{EK1-4}. 
The goal of this paper is to define and study the notion 
of a quantum vertex operator algebra (VOA) in the setting 
of the formal deformation theory and give interesting 
examples of such algebras. Our definition of 
a quantum VOA is based on the ideas 
of the paper \cite{FrR}.

The first chapter of our paper is devoted to the general theory of 
quantum VOAs. For simplicity we consider only bosonic
algebras, but all the definitions and results admit a straightforward 
generalization to the super-case.  

We start with the version of the definition 
of a VOA in which the main axiom is 
the locality (commutativity) axiom. To obtain a quantum deformation 
of this definition, we replace the locality axiom with the 
$\Cal S$-locality axiom, where $\Cal S$ is a shift-invariant unitary 
solution of the quantum Yang-Baxter equation (the other axioms 
are unchanged). We call the obtained structure a braided 
VOA.

However, a braided VOA  
does not necessarily satisfy the associativity 
property, which is one of the main properties of a usual
VOA. More precisely, instead of associativity 
it satisfies a quasi-associativity identity, which differs 
from associativity by insertions of $\Cal S$.  
Therefore we impose an additional axiom, 
which is called the hexagon axiom (by analogy with the theory 
of quasitriangular Hopf algebras). This axiom implies associativity.
We call a braided VOA which satisfies the hexagon
axiom a quantum VOA. Quantum VOAs, 
similarly to the classical ones, possess the operator product expansion. 

Any quantum VOA $V$ has a classical limit $V^0$, 
which is obtained when the quantization parameter $h$ goes to $0$. 
This limit is an ordinary VOA. As usual, 
the classical limit gets a quasiclassical structure, induced by 
the quantum deformation. 
Such a structure is a shift-invariant unitary solution  
$s$ of the classical Yang-Baxter equation with values in 
$\End(V^0)$, which satisfies a quasiclassical version 
of the hexagon relation. It is obtained from $\Cal S$ by 
$s=\frac{d\Cal S}{dh}|_{h=0}$. 

In order to study quasiclassical structures on VOAs, 
we generalize the 
notion of a derivation of a VOA, and introduce the notion
of a pseudoderivation. Pseudoderivations of a VOA 
$V$ form a Lie subalgebra $\text{PDer}(V)$ of 
the Lie algebra $\Hom_k(V,V\o k((z)))$ (here $k$ is the ground field).    

The role of pseudoderivations in the theory of quasiclassical structures 
is that the quasiclassical hexagon 
relation for a tensor $s$ is equivalent to the condition that components 
of $s$ belong $\text{PDer}(V)$; in other words, a quasiclassical structure 
on $V$ is a unitary solution of the classical Yang-Baxter equation
with components in $\text{PDer}(V)$.  
In particular, if $V$ is the affine VOA 
corresponding to a simple Lie algebra $\g$ and any central charge $K$ then 
for any classical r-matrix on $\g$ with spectral parameter one can 
associate a quasiclassical structure on $V$. 

In the second chapter we quantize this quasiclassical 
structure for $\g=sl_N$, when the r-matrix 
on $\g$ is rational, trigonometric, or elliptic. 
This is done using quantum loop groups.
Using the language of the $RTT=TTR$ type relations, we
present the resulting quantum VOA 
by explicit formulas. In this presentation, the simplest vertex 
operator turns out to be the quantum current $\Cal T$ introduced 
by Reshetikhin and Semenov-Tian-Shansky \cite{RS}.
It satisfies the well known $R\Cal TR\Cal T=\Cal TR\Cal TR$-relation 
\cite{RS} which defines the quantum loop algebra.  
These results can be generalized from 
$sl_N$ to any simple Lie algebra.
 
In the third chapter we generalize to the quantum case the well known 
construction of the affine VOA based on considering 
conformal blocks on $\Bbb P^1$ for the Wess-Zunino-Witten model. 
Because of the presence of $\Bbb P^1$, 
this generalization works only for rational r-matrices. 
We show that the obtained quantum VOA 
coincides with the one constructed in the second chapter. 

{\bf Acknowledgements.} The authors are grateful to B.Bakalov, D.Gaitsgory, 
E.Frenkel, I.Frenkel, V.Kac, and B.Lian for useful discussions and references. 
In particular, Prop. 1.9 was pointed out to us by V.Kac. 

\head 1. Quantum vertex operator algebras \endhead

\subhead 1.1. Vertex operator algebras \endsubhead

Vertex operator algebras (VOAs) were introduced by Borcherds \cite{B}. 
There exist many modifications of the definition of a VOA. We are going to 
use the definition from \cite{K}, which we reproduce below. 

Let $k$ be a field of characteristic zero. 
 
\proclaim{Definition} 
A vertex operator algebra (VOA) 
over $k$ is the following data:

1) a $k$-vector space $V$; 

2) a linear map $Y:V\o V\to V((z))$, 
$\v\o \w\to Y(z)(\v\o \w)$;

3) a linear operator $T:V\to V$ (the shift, or Sugawara operator);

4) a vector $\Omega\in V$ (the vacuum vector);

subject to the following axioms:

(A1) locality: 
for any $\v,\w\in V$ there exists $N\ge 0$ such that   
for any $\u\in V$ the  series 
$(z_1-z_2)^NY(z_1)(1\o Y(z_2))(\v\o \w\o \u)$ coincides with the  
series \linebreak $(z_1-z_2)^NY(z_2)(1\o Y(z_1))(\w\o \v\o \u)$. 

(A2) $T\Omega=0$ and $\frac{d}{dz}Y(z)=TY(z)-Y(z)(1\o T)$;

(A3) $Y(z)(\Omega\o \v)=\v$ for all $\v\in V$; 
for any $\v\in V$ the series 
$Y(z)(\v\o \Omega)$ is regular at $z=0$ 
and its value at $z=0$ equals $\v$. 
\endproclaim

{\bf Remark 1.} Note that by axioms (A2) and (A3) 
the operator $T$ is determined by $Y$ and $\Omega$ via 
the formula $T\v=\lim_{z\to 0}\frac{dY(\v,z)\Omega}{dz}$ (cf. \cite{K}).

{\bf Remark 2.} Background material on VOAs can be found 
in several texts: \cite{FLM},\cite{K},\cite{FHL}.
We will use \cite{K} as a main reference for basic facts 
on VOAs. 

In the theory of VOAs, it is convenient to use the following notation. 
For $\v\in V$, denote by $Y(\v,z)$ the element of 
$\End(V)[[z,z^{-1}]]$ defined by the formula
$Y(\v,z)\w=Y(z)(\v\o \w)$. This notation is used in most texts on VOAs and 
will be often used in this paper. 

\subhead 1.2. Example\endsubhead

Now we recall a well-known 
 example of a VOA (see \cite{FZ},\cite{FF2},\cite{L})
which will be central for this paper.  

Let $\g$ be a finite dimensional Lie algebra 
with an invariant inner product $(,)$. 
Let $\hat\g:=\g((t))\oplus kc$ be the corresponding 
centrally extended loop algebra, with the commutation relation 
$$
[a(t),b(t)]=[ab](t)+\text{Res}_{t=0}(da(t),b(t))c, 
$$
and $c$ being the central element. Fix $K\in k$, and 
let $V=U(\hat\g)\o_{U(\g[[t]]\oplus kc)}\chi_K$, 
where $\chi_K$ is the 1-dimensional module over $\g[[t]]\oplus kc$ 
in which $\g[[t]]$ acts by $0$ and $c$ acts by multiplication by $K$. 
This module is called the Weyl module with highest weight $0$.

The space $V$ carries a natural structure 
of a VOA. This structure is defined as follows.

For any $x(z)\in \End(V)[[z,z^{-1}]]$, 
let $x_+(z)$ be the regular part of $x(z)$ with respect to $z$.
Let $x_-(z)=x(z)-x_+(z)$. 
For $a\in \g$, set $a(z)=\sum_{n\in\Z}(a\o t^n)z^{-n-1}\in \hat\g[[z,z^{-1}]]$
 (this expression can be regarded as 
a series with values in $\End(V)$). We have $a(z)=a_+(z)+a_-(z)$. 

Let $\Omega$ be a highest weight vector of $V$. 
The map $Y$ for $V$ is given by the formula
$$
Y(a^1_+(u_1)...a^n_+(u_n)\Omega,z)=
:a^1(z+u_1)...a^n(z+u_n):,\ a^i\in \g\tag 1.1
$$
where $::$ is the normal ordering (cf \cite{K}), 
defined inductively via the formula 
$$
:x^1(z)...x^n(z):=x^1_+(z):x^2(z)...x^n(z):+:x^2(z)...x^n(z):x^1_-(z),
$$
for any $x^i\in \Hom(V,V((z)))$, and 
$a^i(z+u_i):=a^i(z)+u_ia^i(z)'+u_i^2a^i(z)''/2!+...$.  
Formula (1.1) is written in the form of generating functions, so taking 
the coefficients we get a complete description of $Y$. 

The Sugawara operator $T$ is 
the operator $D$ defined by the formula
$$
e^{zD}a^1_+(u_1)...a^n_+(u_n)\Omega=a^1_+(u_1+z)...a^n_+(u_n+z)\Omega\tag 1.2
$$

It is well known that $V$ equipped with $(Y,D,\Omega)$ is a VOA 
(see \cite{L},\cite{K}). 
If $\g$ is abelian and $(,)$ is a nondegenerate form then 
$V$ is called the Heisenberg VOA. 
If $\g$ is a simple Lie algebra and 
$(,)$ the invariant form (normalized so that long roots have squared norm 2)
then $V$ is called the affine VOA. 

{\bf Remark.} 
The VOA corresponding to a reductive Lie algebra $\g$  
was introduced in \cite{FZ,FF2}, based on the ideas 
coming from conformal field theory. Later this definition  
was generalized to any Lie algebra with an inner product in \cite{L}.

In the following sections, the constructed VOA $V$ will be denoted by 
$V(\g,K)$. 

\subhead 1.3. Braided VOAs\endsubhead

{\bf 1.3.1.} Following \cite{FrR}, one may suggest the following
notion of a braided VOA: 

\proclaim{Definition}  A 
braided VOA over $k[[h]]$ is the following data:

1) a topologically free $k[[h]]$-module $V$; 

2) a linear map $Y:V\o V\to V((z))$, 
$\v\o \w\to Y(z)(\v\o \w)$ (here $V((z)):=\{\sum v_nz^n, v_n\to 0\text{ as }
n\to -\infty\}$);

3) a linear operator $T:V\to V$ (the Sugawara operator);

4) a vector $\Omega\in V$ (the vacuum vector);

5) a linear map $\Cal S:V\o V\to V\o V\o k((z))$, 
such that $\Cal S=1+O(h)$ 
which satisfies the shift condition $[T\o 1,\Cal S(z)]=-\frac{d\Cal S}{dz}$, 
the quantum Yang-Baxter equation 
$$
\Cal S^{12}(z)
\Cal S^{13}(z+w)\Cal S^{23}(w)=
\Cal S^{23}(w)\Cal S^{13}(z+w)\Cal S^{12}(z)\tag 1.3
$$
and the unitarity condition 
$$
\Cal S^{21}(z)=\Cal S^{-1}(-z);\tag 1.4
$$

subject to the following axioms:

(QA1) $\Cal S$-locality: 
for any $\v,\w\in V$ and any natural number $M$
there exists $N\ge 0$ such that   
for any $\u\in V$ the  series 
$(z_1-z_2)^NY(z_1)(1\o Y(z_2))(\Cal S(z_1-z_2)(\v\o \w)\o \u)$ 
coincides with the  
series $(z_1-z_2)^NY(z_2)(1\o Y(z_1))(\w\o \v\o \u)$ modulo $h^M$. 

(QA2) $T\Omega=0$ and $\frac{d}{dz}Y(z)=TY(z)-Y(z)(1\o T)$;

(QA3) $Y(z)(\Omega\o \v)=\v$ for all $\v\in V$; 
for any $\v\in V$ the series 
$Y(z)(\v\o \Omega)$ is regular at $z=0$ 
and its value at $z=0$ equals $\v$. 
\endproclaim

{\bf Remark 1.} Here and below, tensor products are understood 
in the h-adically complete sense. 

{\bf Remark 2.} Note that for an infinite dimensional vector space $W$,  
the spaces $W((z))$ and $W\o k((z))$ do not coincide: the second one is 
a proper subspace in the first one. It will be important 
for us to distinguish between them. 

{\bf Remark 3.} Observe that the axioms (QA2),(QA3) coincide with 
(A2),(A3). In particular, 
the operator $T$ is determined by $Y$ and $\Omega$ via 
the same formula as in the classical case. 

\vskip .05in

From the definition of a braided VOA 
it is clear that the reduction $V/hV$ of a braided vertex 
operator algebra $V$ modulo $h$ is a usual vertex operator algebra.
This vertex operator algebra is called the classical limit of $V$. 

{\bf 1.3.2.}
Now we prove a quasi-associativity identity for braided VOAs. 

\proclaim{Proposition 1.1} The map $Y$ satisfies the quasi-associativity 
relation
$$
Y(z)(1\o Y(w))\Cal S^{23}(w)\Cal S^{13}(z)=
Y(w)\Cal S(w)(Y(z-w)\o 1).\tag 1.5
$$
This identity should be understood as follows: on the left 
hand side, $Y(z)$ should be decomposed as $\sum Y^{(m)}(z-w)w^m/m!$, 
and similarly for $\Cal S(z)$, 
and then the equality holds in 
$\Hom(V,V\o V\o V[[w,w^{-1},z-w,(z-w)^{-1}]])$
modulo any given power of $h$ after 
multiplication by a suitable power of $z$.  
\endproclaim

Before proving the Proposition, let us prove the following Lemma. 

\proclaim{Lemma 1.2} For any $\v,\w\in V$, one has 
$$
Y(z)\Cal S(z)(\v\o \w)=e^{zT}Y(-z)(\w\o \v)\tag 1.6
$$
in $V((z))$. 
\endproclaim

\demo{Proof} By axiom (QA2), 
$e^{(z+w)T}Y(-z)(\w\o \v)=Y(w)(\w\o e^{(z+w)T}\v)$ 
(i.e. the inner products of both sides with $f\in V^*$ coincide 
in $k((z))[[z+w,h]]$). 
From axioms (QA2-QA3), we have 
$e^{(z+w)T}\v=Y(z+w)(\v\o \Omega)$. 
Therefore, 
$$
e^{(z+w)T}Y(-z)(\w\o \v)=Y(w)(1\o Y(z+w))(\w\o \v\o \Omega)
$$
(i.e. the inner products of both sides with $f\in V^*$ coincide 
in $k((z))((z+w))[[h]]$). 

Therefore, by $\Cal S$-locality (axiom (QA1)), 
$$
e^{(z+w)T}Y(-z)(\w\o \v)=Y(z+w)(1\o Y(w))(\Cal S(z)(\v\o \w)\o \Omega),
\tag 1.7
$$
i.e. the left and right hand sides
 coincide modulo any given power 
of $h$ after multiplication 
by a suitable power of $z$. 
Taking $w=0$, by axiom (QA3) we get 
$$
e^{zT}Y(-z)(\w\o \v)=Y(z)\Cal S(z)(\v\o \w).
$$
The Lemma is proved. 
$\square$\enddemo

\demo{Proof of Proposition 1.1} 
By Lemma 1.2 we have: 
$$
Y(z)(1\o Y(w))(\Cal S^{23}(w)\Cal S^{13}(z)(\v\o \w\o \u))=
Y(z)(1\o e^{wT}Y(-w))(\Cal S(z)(\v\o \u)\o \w).\tag 1.8
$$
Using (QA2), we get
$$
Y(z)(1\o Y(w))(\Cal S^{23}(w)\Cal S^{13}(z)(\v\o \w\o \u))=
e^{wT}Y(z-w)(1\o Y(-w))
(\Cal S(z)(\v\o 
\u)\o \w).\tag 1.9
$$
On the other hand, consider 
$Y(w)\Cal S(w)(Y(z-w)\o 1)(\v\o \w\o \u)$. 
By Lemma 1.2, we have 
$$
Y(w)\Cal S(w)(Y(z-w)\o 1)(\v\o \w\o \u)=e^{wT}Y(-w)(1\o Y(z-w))
(\u\o \v\o \w).\tag 1.10
$$
The right hand sides of (1.9) and (1.10) coincide by (QA1)
(more precisely, they coincide modulo any given power 
of $h$ after multiplication by a suitable power of $z$). 
Therefore, so do the left hand sides. The Proposition is proved. 
$\square$\enddemo

{\bf Remark.} In the proofs of Lemma 1.2 and Proposition 1.1, 
we followed the logic of the proofs of their classical analogs 
(\cite{FLM}, page 471; \cite{K}), 
in particular the exposition of \cite{FF1}.

\vskip .05in

{\bf 1.3.3.} Recall (\cite{L},\cite{K}) that a derivation of a VOA $V$ is
 a linear map $X:V\to V$ such that $XY(z)=Y(z)(X\o 1+1\o X)$. 
If $X$ is a derivation of $V$ then 
$[T,X]=0$ and $X\Omega=0$. 

\proclaim{Definition} We say that 
a linear map $X:V\to V$ is 
a derivation of a braided VOA $V$ if $XY(z)=Y(z)(X\o 1+1\o X)$ 
and $[X\o 1+1\o X,\Cal S(z)]=0$. 
\endproclaim

It is shown as in the classical case that 
if $X$ is a derivation of $V$ then 
$[T,X]=0$ and $X\Omega=0$. 

Lemma 1.2 implies the following corollary:

\proclaim{Corollary 1.3} In any braided VOA $V$, 
$Y(z)(T\o 1)=\frac{dY(z)}{dz}$; in particular, $T$ is a 
derivation of $V$. 
\endproclaim

\demo{Proof} Indeed, by Lemma 1.2 and (QA2), 
$$
\gather
Y(z)\Cal S(z)(e^{uT}\o 1)=e^{zT}Y(-z)(1\o e^{uT})\sigma=
e^{(z+u)T}Y(-z-u)\sigma=Y(z+u)\Cal S(z+u)=\\
Y(z+u)(e^{-uT}\o 1)\Cal S(z)(e^{uT}\o 1),\endgather
$$
where $\sigma$ is the permutation of two components. 
This implies that $Y(z)(e^{uT}\o 1)=Y(z+u)$, i.e. $Y(z)(T\o 1)=Y'(z)$. 
So, by (QA2) and the property $[T\o 1,\Cal S]=-\Cal S'$, 
we get that $T$ is a derivation of $V$. 
$\square$\enddemo

\subhead 1.4. Quantum VOAs\endsubhead

{\bf 1.4.1.}
In the last section we saw that in general braided VOAs satisfy a 
quasi-associativity identity which is different from the classical 
associativity identity. However, from many points of view (physical 
and mathematical), it is natural to consider 
associative algebras only. 
Therefore, we will make the following definition, which is motivated by 
Proposition 1.4 below.  

\proclaim{Definition} A braided VOA $V$ is said to be a quantum VOA if 
the following hexagon relation is satisfied:

(QA4)
$$
\Cal S(w)(Y(u)\o 1)=(Y(u)\o 1)\Cal S^{23}(w)\Cal S^{13}(u+w)\tag 1.11
$$
(that is, matrix elements of the two sides of (1.11) coincide in 
$k((w))((u))[[h]]$). 
\endproclaim

\proclaim{Proposition 1.4} In a quantum VOA the following associativity 
axiom holds:
$$
Y(z)(1\o Y(w))=
Y(w)(Y(z-w)\o 1).\tag 1.12
$$
This equality is understood as in Proposition 1.1.
\endproclaim

\demo{Proof} By the hexagon relation, 
$\Cal S(w)(Y(z-w)\o 1)=(Y(z-w)\o 1)\Cal S^{23}(w)\Cal S^{13}(z)$. 
Therefore, we get from Proposition 1.1
$$
Y(z)(1\o Y(w))\Cal S^{23}(w)\Cal S^{13}(z)=
Y(w)(Y(z-w)\o 1)\Cal S^{23}(w)\Cal S^{13}(z).\tag 1.13
$$
Since $\Cal S$ is invertible, the Proposition follows. 
$\square$\enddemo

\proclaim{Corollary 1.5} For any $\v,\w$ in a quantum VOA $V$ 
and any integer $n$ there exists a unique 
$\y_n\in V$ such that $\text{Res}_{z=w}(z-w)^nY(z)(1\o Y(w))
(\v\o \w \o\u)
=Y(w)(\y_n\o \u)$. 
\endproclaim

\demo{Proof} Indeed, by Proposition 1.4, 
the vector $\y_n=\Res_{z=0}z^nY(z)(\v\o\w)$ satisfies the required conditions.
Suhc a vector is obviously unique.  
$\square$\enddemo

Corollary 1.5 implies the existence of the operator product expansion
for quantum VOA:  
$$
Y(\v,z)Y(\w,w)=\sum_n(z-w)^{-n-1}Y(\y_n,w).\tag 1.14
$$

{\bf Remark.} Proposition 1.4. 
and corollary 1.5 show that quantum VOAs are a special case 
of field algebras over $k[[h]]$
(cf. \cite{K}). More specifically, a quantum VOA is a field 
algebra with an R-matrix. 

{\bf 1.4.2.} 
\proclaim{Proposition 1.6} In a quantum VOA, $\Cal S(z)(\Omega\o \v)
=\Omega\o \v$.
\endproclaim

\demo{Proof} Represent $V$ in the form $V^0[[h]]$, where $V^0$ 
is a vector space, 
so that 
$\Omega$ is h-independent (this can be done without loss of generality). 
Then $\Cal S(z)=1+\sum_{i\ge 1}h^is_i(z)$, so it is enough 
to show that $s_i(z)(\Omega\o \v)=0$. We show it by induction in $i$. 
The base and step of induction are justified by the same argument. 
Suppose for $j<i$ the statement is known (it is definitely so for $i=1$), 
and let us prove it for $j=i$. Let $Y=Y_0+hY_1+...+h^nY_n+...$
Computing the coefficient of 
$h^i$ in (QA4), we get
$$
\gather
s_i(w)(Y_0(u)\o 1)+\sum_{j<i}s_j(w)(Y_{i-j}(u)\o 1)=\\
(Y_0(u)\o 1)(s_i^{23}(w)+s_i^{13}(w+u))+\sum_{j_1+j_2<i}
(Y_{i-j_1-j_2}(u)\o 1)s_{j_1}^{23}(w)
s_{j_2}^{13}(w+u).\endgather
$$
When we apply both sides of this equation to $\Omega\o \w\o \v$, 
all terms but one drop out, and we get
$$
(Y_0(u)\o 1)s_i^{13}(w+u)(\Omega\o \w\o \v)=0
$$
Let $f\in (V^0)^*$, and $\sigma(w)=(1\o f)(s_i(w)(1\o \v))$. We have 
$Y_0(\sigma(w+u)\w,u)=0$. 
But it is easy to see that if $\v(z)\in V^0\o k((z))$ and 
$Y(\v(z),z)=0$ then $\v(z)=0$: 
indeed, if $\v(z)=\v_nz^n+$ higher degree terms, then 
$\lim_{z\to 0}z^{-n}Y(\v(z),z)=\v_n$. 
This implies that $\sigma(w)\w=0$ for all $\w$. 
The proposition is proved.
$\square$\enddemo

\subhead 1.5. Quasiclassical VOAs\endsubhead

Now consider the quasiclassical object corresponding 
to the notion of a quantum VOA. 

\proclaim{Definition} Let $V$ be a vertex operator algebra.
A classical r-matrix on $V$ 
is a linear map $s(z):V\o V\to V\o V\o k((z))$, satisfying 
the following conditions:

(i) the classical Yang-Baxter equation with spectral parameter:
$$
[s^{12}(z_1-z_2),s^{13}(z_1-z_3)]+
[s^{12}(z_1-z_2),s^{23}(z_2-z_3)]+
[s^{13}(z_1-z_3),s^{23}(z_2-z_3)]=0;\tag 1.15
$$

(ii) the unitarity condition 
$$
s^{21}(-z)=-s(z);\tag 1.16
$$

(iii) the shift condition $[T\o 1,s(z)]=-\frac{ds}{dz}$; and 

(iv) the hexagon relation 
$$
s(w)(Y(u)\o 1)=(Y(u)\o 1)(s^{23}(w)+s^{13}(u+w)).\tag 1.17
$$

A quasiclassical VOA is a VOA equipped with a classical r-matrix. 
\endproclaim

Let $V$ be a quantum VOA. Then the VOA $V^0=V/hV$ 
has a natural classical $r$-matrix  $s(z)$ 
defined by $\Cal S(z)=1+hs(z)+O(h^2)$. The corresponding quasiclassical
VOA $V^0$ is said to be the quasiclassical limit of the quantum VOA $V$, and 
$V$ is said to be a quantization of $(V^0,s)$. 

It is an interesting question whether any quasiclassical VOA can be quantized. 
In the next section, we show that it is so in the trivial case 
when the Sugawara element vanishes.

\subhead 1.6. A toy example \endsubhead

{\bf 1.6.1.} A useful toy example of the above theory is the case when $T=0$ 
(``topological field theory''). 
It is obvious from the definition of a VOA that a VOA $V$ 
with $T=0$ is the same 
thing as a 
commutative associative algebra with unit: the map $Y$ does not 
depend on $z$ and defines a commutative associative product. 

Further, a classical r-matrix on a VOA $V$ with $T=0$ is 
simply a
skewsymmetric solution $s\in \End(V\o V)$ 
of the classical Yang-Baxter equation 
(without spectral parameter) such that for any $\v\in V$, $f\in V^*$ 
the element $(1\o f)(s(1\o \v))\in \End(V)$ is a derivation of $V$ as an 
algebra. If $V$ is finitely generated, then this implies that 
$r\in \Lambda^2\text{Der}(V)$, where $\text{Der}(V)$ is the Lie algebra 
of derivations of $V$. 

Finally, a quantum VOA with $T=0$ is an associative unital algebra $V$ 
equipped with an element $\Cal S=1+hs+O(h^2)\in \End(V\o V)[[h]]$ 
which satisfies the quantum Yang-Baxter equation and the unitarity 
condition, such that $(Y\o 1)\Cal S^{23}\Cal S^{13}=
\Cal S(Y\o 1)$ 

{\bf Remark.}
Notice that a braided VOA 
with $T=0$ which does not satisfy (QA4) need not be associative, 
and thus is not a reasonable object. This is one of our motivations to 
impose axiom (QA4).  

{\bf 1.6.2.}
\proclaim{Proposition 1.7} Let $V^0$ be a finitely generated
commutative associative algebra. 
Then any structure $s$ of a quasiclassical VOA on $V$ can be quantized.
\endproclaim  

\demo{Proof} Let $\g$ denote the Lie algebra of derivations of $V$. 
Let $U_h(\g)$ be the quantization of the 
triangular Lie bialgebra $(\g,s)$ as constructed in \cite{EK1}. 
This quantization is equal to $U(\g)[[h]]$ as an algebra, 
with the coproduct given by $\Delta(x)=J^{-1}\Delta_0(x)J$, 
where $J=1+hs/2+O(h^2)$ is a multiplicative 2-cocycle, i.e. 
$$
(\Delta_0\o 1)(J)(J\o 1)=(1\o \Delta_0)(J)(1\o J),\
(\e\o 1)(J)=(1\o \e)(J)=1,\tag 1.18
$$
where $\Delta_0,\e$ are the coproduct and counit in $U(\g)$.  
Define $V=V^0[[h]]$, 
and let $\Cal S=J_{21}^{-1}J$ be the universal 
R-matrix of $U_h(\g)$. Set $Y=Y^0J_{21}$, 
where $Y^0$ is the product in $V^0$. 
It is easy to check that this data  
 satisfies the axioms of a braided VOA (with $T=0$). 
To check that $Y$ satisfies (QA4), 
observe that $Y:V\o V\to V$ is a morphism in the category of 
$U_h(\g)^{op}$-modules. Since $(\Delta^{op}\o 1)(\Cal S)=\Cal S^{23}
\Cal S^{13}$, 
we see that (QA4) holds. Thus, $V$ is a quantum VOA. It is easy to see that 
$V$ is a quantization of $V^0$.  
$\square$\enddemo

{\bf 1.6.3.}
This example of quantization has the following 
straightforward generalization. Suppose 
$V$ is any VOA (with $T$ not necessarily $0$), 
$\g$ the Lie algebra of derivations of $V$, and 
$s\in \Lambda^2\g$ is a skewsymmetric 
classical r-matrix.  
Then $s$ defines a structure of a quasiclassical 
VOA on $V$. Moreover, any structure of a 
quasiclassical VOA which is independent on $z$ 
is of this form if $V$ is finitely generated
(for the definition of a finitely generated VOA, see section 1.7). 
This structure can be quantized as described 
in the proof of Proposition 1.6. 

\subhead 1.7. Pseudoderivations of VOAs and quasiclassical 
VOAs\endsubhead
  
{\bf 1.7.1.}
Now consider a classical r-matrix $s(z)$ on a VOA $V$ which is not necessarily 
constant. 

Let $\v\in V$, $f\in V^*$, and $a(u)=(1\o f)(s(u)(1\o \v))\in \Hom_k(V,V\o 
k((u)))$. 
By definition, 
$$
[T,a(u)]=-\frac{da}{du},\tag 1.19
$$ 
Also, we have
$$
[a(u),Y(\w,z)]=Y(a(z+u)\w,z),\tag 1.20
$$
where $a(z+u)$ is expanded in the Taylor series in $z$.

This gives rise to the following definition: 

\proclaim{Definition} An element $a\in \Hom(V,V\o k((z)))$ is said to be 
a pseudoderivation of a VOA $V$ if it satisfies (1.19) and (1.20). 
\endproclaim

Thus, the hexagon axiom for a classical $r$-matrix $s$ on a VOA $V$ 
is equivalent to saying 
that for any $\v\in V,f\in V^*$ the element $(1\o f)(s(1\o \v))$
is a pseudoderivation of $V$. 

It is clear that pseudoderivations of $V$ form a Lie algebra $\text{PDer}(V)$, 
which contains the Lie algebra $\text{Der}(V)$ of derivations. 
More specifically, the operator 
$\text{ad}(T)=-\frac{d}{dz}$ 
is a derivation of $\text{PDer}(V)$, and the space of its invariants 
(i.e. constant elements) is $\text{Der}(V)$.

{\bf 1.7.2.} 
Let us consider some properties of pseudoderivations. 

Recall first that a VOA $V$ is said to be generated 
by a set $S\subset V$ if any vector in $V$ can be obtained from $S$ 
by iterating the following operations: 

1) linear combination;

2) the operator $T$;

3) the operation $\v\circ_n\w:= \Res_{z=0}z^{-n-1}Y(\v,z)w$ for any 
integer $n$. 

If $S$ is finite, $V$ is said to be finitely generated. 

\proclaim{Proposition 1.8} (i) For any pseudoderivation $a(z)$ of a 
VOA $V$, 
one has $a\Omega=0$.

(ii) Let $V$ be a VOA generated by $S$. 
Then a pseudoderivation of $V$ is completely determined by its action 
on $S$. 
\endproclaim

\demo{Proof} (i) Recall from 
the proof of Proposition 1.6 that if $\v(z)\in V\o k((z))$ 
and $Y(\v(z),z)=0$ then $\v=0$. 
By (1.20) we have $Y(a(z+u)\Omega,z)=0$. 
By (1.19), $a(z+u)\Omega\in V\o k((u))[[z]]$. 
Therefore, $a(u)\Omega=0$. 

 (ii) By (1.20), the action of $a$ on $S$ 
completely determines the commutators $[a(z),Y(\v,w)]$ for $\v\in S$. 
Further, it is easy to see that if $[a(z),Y(\v,w)]$ is known for input vectors
of operations 1-3, then its is known for the output vector. 
Indeed, this is obvious for operation 1, follows from $Y(T\v,z)=
Y'(\v,z)$ for operation 2, and follows from associativity for operation 3. 
Thus, $[a(z),Y(\v,w)]$ is known for all $\v\in V$. This means (by (1.20)) 
that $Y(a(z+w)\v,w)$ is known for all $\v$. As in the 
proof of (i), this implies that $a(z+w)\v$ is known 
for all $\v$, or, setting $z=0$, $a(w)\v$ is known for all $\v$. 
$\square$\enddemo

{\bf 1.7.3.} The main source of examples of pseudoderivations is the following 
proposition. 

\proclaim{Proposition 1.9} 
Let $V$ be a VOA, and $Y(\v,z)=\sum Y_n(\v)z^{-n-1}$ for $\v\in V$ 
be a vertex operator. Then for any $\alpha\in k((z))$ 
the element 
$$
X_{\alpha,\v}:=\Res_{w=0}\alpha(z+w)Y(\v,w)dw=
\sum_{m\ge 0}\frac{\alpha^{(m)}(z)Y_m(\v)}{m!}\tag 1.21
$$
is a pseudoderivation of $V$. 
\endproclaim

\demo{Proof} Identity (1.19) 
is obvious. To prove 
identity (1.20), we will use Borcherds' commutator formula (see \cite{K}):
if $x(z)=\sum x_mz^{-m-1}$ is a vertex operator then 
$$
[x_m,Y(\v,z)]=\sum_{j\ge 0}\left(\matrix m\\ j\endmatrix\right)
Y(x_j\v,z)z^{m-j}
$$
Using this formula, we get
$$
\gather
[\sum_{m\ge 0}\frac{\alpha^{(m)}(z)Y_m(\v)}{m!},Y(\w,w)]=
\sum_{m\ge j\ge 0}\left(\matrix m\\ j\endmatrix\right)
\frac{\alpha^{(m)}(z)w^{m-j}Y(Y_j(\v)\w,w)}{m!}=\\
\sum_{m\ge j\ge 0}
\frac{\alpha^{(m)}(z)w^{m-j}Y(Y_j(\v)\w,w)}{(m-j)!j!}=
\sum_{j\ge 0}
\frac{\alpha^{(j)}(z)Y(Y_j(\v)\w,w)}{j!},
\endgather
$$
which implies (1.20). 
$\square$\enddemo

\proclaim{Corollary 1.10} Let $V=V(\g,K)$. Then 
for any $a\in \g$, and 
any $\alpha\in k((z))$ the element 
$$
a\o \alpha(z+t):=\sum_{m\ge 0}\frac{\alpha^{(m)}(z)}{m!}a\o t^m.\tag 1.22
$$
is a pseudoderivation of $V$.
\endproclaim

\demo{Proof} Follows from Proposition 1.9, since $a(z)=Y(a_{-1}\Omega,z)$ is 
a vertex operator. 
$\square$\enddemo

\subhead 1.8. Examples of quasiclassical VOAs \endsubhead

{\bf 1.8.1.}
Now consider an example of a quasiclassical VOA. 
Let $V=V(\g,K)$. By Corollary 1.10, 
the Lie algebra $\g((t))$ is a Lie 
subalgebra of $\text{PDer}(V)$: the embedding 
$\phi:\g((t))\to \text{PDer}(V)$ is given by the formula 
$a\o \alpha(t)\to a\o \alpha(z+t)$. 

Let us look for classical 
r-matrices $s$ whose components 
(i.e. expressions \linebreak 
$(1\o f)(s(1\o \v))$) are in this Lie subalgebra.  
It is easy to show that any such $s$ has the form
$$
s(z)=\rho(u-v+z),\tag 1.23
$$
where $\rho\in \g\o\g((z))$ is an element satisfying the 
classical Yang-Baxter equation and the unitarity condition.
Here $\rho(u-v+z)$ is regarded as an element of 
$\g[u]\o\g[v]((z))$ (i.e. we are expanding 
$\rho(u-v+z)$ in a Taylor series with respect to $u$ and $v$). 

Thus, we conclude that a classical r-matrix on $V$ whose components 
are in $\g((t))$ is the same thing as a classical r-matrix 
with a spectral parameter with values in $\g$. 

{\bf 1.8.2.} 
Consider the special case when $\g$ is abelian. 
Then the quasiclassical 
structure $s$ admits an easy quantization, using the ordinary 
exponential. Indeed, define 
$J(z)=e^{h\rho(u-v+z)/2}$, and $\Cal S(z)=J_{21}^{-1}(-z)J(z)=
e^{h\rho(u-v+z)}$. Define a new vertex operator product on $V$ by 
$\tilde Y(z)=Y(z)J_{21}(-z)$. It is easy to see that 
$(V[[h]],\tilde Y,T,\Omega,\Cal S)$ is a quantum VOA which is a quantization
of the quasiclassical VOA $(V,Y,T,\Omega,s)$. This method actually works 
for any classical r-matrix whose components are in a commutative subalgebra 
of $\g((z))$. 

As an example consider the case when $\g$ is 1-dimensional,
with basis vector $a$ such that $(a,a)=1$. In this case 
$\rho$ is any  scalar-valued odd element of $k((z))$.
Define $\tilde a(z)=\tilde Y(a_{-1}\Omega,z)$. Then 
it is easy to check that 
$$
\tilde a(z)=a(z)-\frac{h}{2}\Res_{u=0}\rho'(z-u)a(u).\tag 1.24
$$
 
However, if $\g$ is non-abelian, quantization is a more delicate matter. 
It is discussed in Chapter 2 in the special case $\g=sl_N$.  

\subhead 1.9. Nondegenerate VOA\endsubhead

{\bf 1.9.1.} 
\proclaim{Definition} A VOA $V$ is said to be nondegenerate 
if the maps 
$$
Z_n=Y(z_1)(1\o Y(z_2))....(1^{\o n-1}\o Y(z_n))(1^{\o n}\o\Omega): 
V^{\o n}\o k(z_1,...,z_n)\to V((z_1))...((z_n))
$$ are injective for all $n$. 
\endproclaim

This notion is useful because of the following proposition. 

\proclaim{Proposition 1.11} Let $(V,Y,T,\Omega, \Cal S)$ be 
a data satisfying the axioms of a braided VOA, except maybe equations 
(1.3) and (1.4). 
Suppose also that $V^0=V/hV$ is a nondegenerate VOA. 
Then: 

(i) equations (1.3) and (1.4) are automatically satisfied;

(ii) axiom (QA4) is equivalent to the associativity identity (1.12).  
\endproclaim

\demo{Proof} (i) Equation (1.3) holds after applying $Z_3$ by (QA1)
(just compare two ways of rewriting the product of three vertex
operators in the reverse order), so it 
must be satisfied, as $Z_3$ is injective. 
Similarly, equation (1.4) is satisfied after applying $Z_2$, so it must 
be satisfied. 

(ii) The associativity relation (1.12) implies (1.13). Combining (1.13) 
with Proposition 1.1, one gets that (QA4) is satisfied after applying 
$Y$. By nondegeneracy, (QA4) must be satisfied. 
$\square$\enddemo

{\bf 1.9.2.} 
\proclaim{Proposition 1.12} 
Let $V=V(\g,K)$. 
If $V$ is an irreducible $\hat\g$-module then 
$V$ is a nondegenerate VOA. 
\endproclaim

\demo{Proof} Let $X\in \text{Ker}Z_n$
be of lowest degree with respect to the natural grading. 
Applying $a_-(u)$ ($a\in \g$) to the identity $Z_nX=0$, we get
$Z_n(a_-^1(u+z_1)+...+a_-^n(u+z_n))X=0$. 
Set $W=(a_-^1(u+z_1)+...+a_-^n(u+z_n))X$. 
The degree of $W$ is not higher than that of $X$, 
and up to lower degree terms we have $W=\sum \frac{a^j_0}{u+z_j}X+..$. 
Consider 
$P=(\frac{d}{du})^n[(u+z_1)...(u+z_n)W]$. Then $P$ is of lower 
degree than $X$, and $Z_nP=0$, so $P=0$ by the assumption of minimality 
of degree. 
This implies that $W=\frac{Q(u)}{(u+z_1)...(u+z_n)}$, where $Q$ 
is a polynomial of $u$ of degree $\le n-1$. This implies 
that coefficients of $\frac{1}{(u+z_j)^{i+1}}$ for $i\ge 1$ 
in the partial fraction expansion of $W$ must be zero. 
Thus, $a_i^jX=0$ for $i\ge 1$ and all $j$. Since $V$ is an irreducible 
module over $\hat\g$, we get that 
$X=f\Omega^{\o n}$. This implies the proposition, as
for such $X$ it is clear that $Z_nX=0$ implies $f=0$. 
$\square$\enddemo

\proclaim{Proposition 1.13} Let the inner product $(,)$ on $\g$ be 
nondegenerate, and $V=V(\g,K)$. Then $V$ is nondegenerate
for generic $K$.  
\endproclaim

\demo{Proof} By Proposition 1.12, it is enough to show that
 $V$ is an irreducible $\hat\g$-module 
for a generic $K$. 

As $K\to\infty$, the Lie algebra $\hat\g$ 
turns into the Heisenberg algebra, and $V$ tends to the Fock module, 
which is irreducible. This implies the irreducibility of $V$ for generic $K$.
$\square$\enddemo

\head 2. The quantum affine VOA\endhead

\subhead 2.1. The construction\endsubhead

In this section we will give a nontrivial example of a quantum vertex 
operator algebra. We will use definitions and notation from \cite{EK4}

{\bf 2.1.1.}
Let $R(u)$ be a rational, trigonometric, or elliptic R-matrix
from \cite{EK4}, which has been normalized 
to satisfy the crossing symmetry condition as in Section 1.4 of \cite{EK4}
(in \cite{EK4} it was denoted by $\bar R$). We recall from \cite{EK4} that the 
quantum function algebra $F_0(R)$ is generated by 
 the entries
of the coefficients of 
the formal series $T(u)\in \End(k^N)\o F(R)[[u]]$,
 $T(u)=T_0+T_1u+...$, 
and the defining relations are
$$
\gather
R^{12}(u-v)T^{13}(u)T^{23}(v)=T^{23}(v)T^{13}(u)R^{12}(u-v),\\
\text{qdet}_R(T(u))=1\tag 2.1\endgather
$$

{\bf Remark.} The second relation in (2.1), as well as similar relations 
below, should be understood as in Proposition 1.1. 

This algebra is a (topological) Hopf algebra, with the coproduct, 
counit, and antipode defined by
$$
\Delta(T(u))=T^{12}(u)T^{13}(u),\ \e(T(u))=
1,\ S(T(u))=T^{-1}(u).\tag 2.2
$$

Let $U_0(R)$ be the QUE algebra obtained by extending the quantum function 
algebra $F_0(R)$. It is generated by $t(u)=\frac{T(u)-1}{h}$ with the 
relations obtained from (2.1) by dividing by the smallest power of $h$ which 
occurs after rewriting the equations in terms of $t(u)$:
$$
\gather
[t^{13}(u),t^{23}(v)]=\\
[r_*^{12}(u-v),t^{13}(u)+t^{23}(v)]
+h(r_*^{12}(u-v)t^{13}(u)t^{23}(v)-t^{23}(v)t^{13}(u)r_*^{12}(u-v)),\\
h^{-1}(\text{qdet}_R(1+ht(u))-1)=0,\endgather
$$
where 
$r_*=h^{-1}(1-R)$.

{\bf 2.1.2.}
Set $V=U_0(R)$. Let $K\in k$. 
In the following, we will introduce the structure of a quantum VOA on $V$, 
such that $V/hV$ is the affine VOA for $\g=sl_N$ at level $K$.  

First of all, let the vacuum vector $\Omega$ be the unit $1$ of the algebra 
$U_0(R)$. Next, define the Sugawara element $D$ on $V$ by
$$
e^{zD}T(u_1)...T(u_n)\Omega=T(u_1+z)...T(u_n+z)\Omega.\tag 2.3
$$
(the series $T(u)$ is an analogue of 
$a_+(u)$, and equation (2.3) is an analogue of (1.2) in the classical case). 
Using defining relations (2.1) for $U_0(R)$, it is easy to see that 
$D$ is well defined. 

Now comes the main challenge -- defining $Y$ and $\Cal S$.
Before we define them, we will introduce the Laurent series $T^*(u)$ of
``annihilation operators'' (as opposed to ``creation operators'' 
$T(u)$). 

\proclaim{Lemma 2.1} 
There exists a unique operator series
 $T^*(u)\in \End(V)[[u,u^{-1}]]$ such that 
$$
\gather
T^{*0,n+1}(u)T^{1,n+1}(v_1)...T^{n,n+1}(v_n)\Omega=
R^{10}(u-v_1-Kh/2)...R^{n0}(u-v_n-Kh/2)\times \\
T^{1,n+1}(v_1)...T^{n,n+1}(v_n)R^{n0}(u-v_n+Kh/2)^{-1}...
 R^{10}(u-v_1+Kh/2)^{-1}\Omega\tag 2.4\endgather
$$
(here $ R$ is 
regarded as a matrix with entries 
in $k((u))[[v,h]]$). 
\endproclaim

\demo{Proof} Formula (2.4) defines $T^*$ on the free algebra generated by 
$t(u)$. So we need to show that $T^*$ descends to the quotient 
by the defining relations, i.e. that it maps the ideal of relations to 
itself. So we have to check that the expression
$$
\gather
 T^{*0,n+1}(u)\biggl(T^{1,n+1}(v_1)...R^{ii+1}(v_i-v_{i+1})
T^{i,n+1}(v_i)T^{i+1,n+1}(v_{i+1})...T^{n,n+1}(v_n)\Omega-\\
T^{1,n+1}(v_1)...
T^{i+1,n+1}(v_{i+1})T^{i,n+1}(v_i)R^{ii+1}(v_i-v_{i+1})
...T^{n,n+1}(v_n)\Omega\biggr)\endgather
$$
belongs to the ideal of relations. 

This statement can be checked by a direct computation. 
For simplicity we check only the case $n=2$; the general case 
is completely analogous. We have
$$
\gather
 T^{*0,3}(u)\biggl(R^{12}(v_1-v_2)
T^{13}(v_1)T^{23}(v_2)-
T^{23}(v_2)T^{13}(v_1)R^{12}(v_1-v_2)\biggr)\Omega=\\
R^{12}(v_1-v_2)R^{10}(u-v_1-Kh/2)R^{20}(u-v_2-Kh/2)\times \\
T^{13}(v_1)T^{23}(v_2)R^{20}(u-v_2+Kh/2)^{-1}R^{10}(u-v_1+Kh/2)^{-1}\Omega-\\
R^{20}(u-v_2-Kh/2)R^{10}(u-v_1-Kh/2)\times \\
T^{23}(v_2)T^{13}(v_1)
R^{10}(u-v_1+Kh/2)^{-1}R^{20}(u-v_2+Kh/2)^{-1}R^{12}(v_1-v_2)\Omega.
\endgather
$$
Using the Yang-Baxter equation for $R$, 
we find that the last expression equals to 
$$
\gather
R^{20}(u-v_2-Kh/2)R^{10}(u-v_1-Kh/2)\biggl(R^{12}(v_1-v_2)
T^{13}(v_1)T^{23}(v_2)-\\
T^{23}(v_2)T^{13}(v_1)
R^{12}(v_1-v_2)\biggr)R^{20}(u-v_2+Kh/2)^{-1}R^{10}(u-v_1+Kh/2)^{-1}\Omega.
\endgather
$$
The expression on the right hand side belongs to the ideal of relations. 

The fact that the determinant identity is preserved by $T^*(u)$ is also 
checked directly. The Lemma is proved. 
$\square$\enddemo

{\bf Remark.} 
In general, $T^*(u)$ is a series which is infinite in both directions, 
but for any $\v\in V$ the series $T^*(u)\v$ is finite in the negative 
direction modulo any given power of $h$. 
Also, in the rational case $T^*(u)$ has only nonpositive powers of $u$ 
and has constant term $1$.   

\proclaim{Proposition 2.2}
The operators $T(u),T^*(u)$ satisfy the following commutation relations:
$$
R^{12}(u-v)T^{*13}(u)T^{*23}(v)=T^{*23}(v)T^{*13}(u)R^{12}(u-v),\tag 2.5
$$
and 
$$
 R^{12}(u-v-Kh/2)T^{13}(u)T^{*23}(v)=T^{*23}(v)T^{13}(u)
 R^{12}(u-v+Kh/2).
\tag 2.6
$$
\endproclaim

\demo{Proof} Straightforward.$\square$\enddemo

{\bf Remark.} 
This proposition shows that $V$ admits a natural 
action of the centrally extended 
quantum double of $U_0(R)$.  

Now define some convenient notation.

For $m,n>0$ and $u=(u_1,...,u_n),v=(v_1,...,v_m)$, 
where $u_i,v_j$ are variables, define 
$$
 R_{nm}(u|v|z)=
\prod_{i=n+m}^{n+1}\prod_{j=1}^n R^{ji}(u_j-v_i+z)\tag 2.7
$$
(note the order of the factors!)
Also, define 
$$
\gather
T_n(u|z)=
T^{1,n+1}(u_1+z)...T^{n,n+1}(u_n+z),\\  
T^*_n(u|z)=
T^{*1,n+1}(u_1+z)...T^{*n,n+1}(u_n+z).
\tag 2.8\endgather
$$ 
These expressions satisfy the following 
commutation relations:
$$
\gather
 R^{12}_{nm}(u|v|z-w)
T^{13}_{n}(u|z)
T^{23}_{m}(v|w) =
T^{23}_{m}(v|w) 
T^{13}_{n}(u|z)
 R^{12}_{nm}(u|v|z-w),\\
 R^{12}_{nm}(u|v|z-w)
T^{*13}_{n}(u|z)
T^{*23}_{m}(v|w) =
T^{*23}_{m}(v|w) 
T^{*13}_{n}(u|z)
 R^{12}_{nm}(u|v|z-w),\\
 R^{12}_{nm}(u|v|z-w-Kh/2)
T^{13}_{n}(u|z)
T^{*23}_{m}(v|w) =
T^{*23}_{m}(v|w) 
T^{13}_{n}(u|z)
 R^{12}_{nm}(u|v|z-w+Kh/2),\tag 2.9
\endgather
$$ 

Now define $Y:V\o V\to V((z))$ by the formula 
$$
Y(T_n(u|0)\Omega,z)=
T_n(u|z)T^*_n(u|z+Kh/2)^{-1}.
\tag 2.10
$$
(note that the fact that the annihilation operators $T^*$
are positioned on the right from $T$ is analogous to the normal ordering in 
the classical formula (1.1)). 
Similarly to Lemma 2.1, 
it is easy to check using the defining relations (2.1) 
that $Y$ is well defined.

Now define $\Cal S:V\o V\to V\o V\o k((z))$
 by the formula
$$
\gather
\Cal S^{34}(z)
\left( R^{12}_{nm}(u|v|z)^{-1}
T^{24}_m(v|0)
 R^{12}_{nm}(u|v|z-Kh)
T^{13}_n(u|0)(\Omega\o \Omega)\right)=\\
T^{13}_n(u|0)
 R^{12}_{nm}(u|v|z+Kh)^{-1}
T^{24}_m(v|0)
 R_{nm}^{12}(u|v|z).\tag 2.11
\endgather
$$
As before, it is easy to check directly that $\Cal S$ is well defined. 

{\bf 2.1.3.} 
\proclaim{Theorem 2.3}
The space $V$ equipped with $(Y,D,\Omega,\Cal S)$ is a quantum VOA. 
\endproclaim

\demo{Proof} Axioms (QA2),(QA3) 
and the commutation relation between 
$D$ and $\Cal S$ are obvious from the definitions. 
So it remains to prove (QA1) and (QA4). 

Axiom (QA1) follows from commutation relations (2.9). 
The Yang-Baxter equation and unitarity condition for $\Cal S$ 
and axiom (QA4) are checked by a 
direct computation.
$\square$\enddemo

{\bf 2.1.4.} 
Note that formula (2.10) implies $Y(T(u),z)=\Cal T(u+z)$, where 
$$
\Cal T(z)=T(z)T^*(z+Kh/2)^{-1}.
$$
The series $\Cal T(z)$ 
has appeared in the literature before 
and is called the quantum current \cite{RS}. It satisfies
the following commutation relation
(formula (2) in \cite{RS}), which is easy to check using (2.9):
$$
\gather
\Cal T^{13}(z_1)
 R^{12}(z_1-z_2+Kh)^{-1}\Cal T^{23}(z_2)
 R^{12}(z_1-z_2)=\\
 R^{12}(z_1-z_2)^{-1}
\Cal T^{23}(z_2) R^{12}(z_1-z_2-Kh)\Cal T^{13}(z_1).\tag 2.12
\endgather
$$
The (suitably completed) 
algebra generated by $\Cal T$ with this relation 
and the unimodularity relation is the 
centrally extended double of $F_0(R)$ modulo the relation 
$c=K$.

\subhead 2.2. The quasiclassical limit of $V$\endsubhead

\proclaim{Proposition 2.4} The VOA $V^0=V/hV$ is isomorphic to 
$V(sl_N,K)$. 
The quasiclassical structure $s$ on $V^0$ defined by $V$ 
is given by formula (1.23), where $\rho=\frac{d R}{dh}|_{h=0}$. 
\endproclaim

\demo{Proof} 
Let $ R=1-hr+O(h^2)$, $T=1+ht+O(h^2)$, 
$T^*=1+ht^*+O(h^2)$. Then the commutation relations for $T,T^*$ 
imply the following commutation relation for $t,t^*$
on $V^0$: 
$$
\gather
[t^{13}(u),t^{23}(v)]=[r^{12}(u-v),t^{13}(u)+t^{23}(v)],\\
[t^{*13}(u),t^{*23}(v)]=[r^{12}(u-v),t^{*13}(u)+t^{*23}(v)],\\
[t^{13}(u),t^{*23}(v)]=[r^{12}(u-v),t^{13}(u)+t^{*23}(v)]+Kr'(u-v)^{12},
\tag 2.13\endgather
$$

Let $\g=sl_N$, and $\g(r)$ be the Lie algebra associated to 
the classical $r$-matrix $r$ in Chapter 2 of \cite{EK3}.
Recall from \cite{EK3} that we can include $\g(r)$ into  
a Manin triple $(\g((x)),\g(r),\g[[x]])$. 
Equations (2.13) imply that $V^0=U(\g(r))$, and the operators $t,t_*$ 
define an action of the affine Lie algebra 
$\hat\g$ on $V^0$ with central charge $K$
(more precisely, 
$t$ defines the action of $\g(r)$ and $t^*$ the action of $\g[[x]]$). 
It is easy to see that as a $\hat\g$-module, $V^0$ is isomorphic 
to the Weyl module with highest weight $0$ considered in Section 1.2, 
and the isomorphism is unique up to scaling. After identifying 
$V^0$ with the Weyl module, the formulas for $t$ and $t^*$  
are given by
$$
t(u)=r^{21}(x-u),\ t^*(u)=-r^{21}(x-u),\tag 2.14
$$
where in the first formula the expansion is in the nonnegative powers of $u$ 
and in the second one in the nonnegative powers of $x$
(in both cases the right hand side is regarded as 
an element of $\g\o\hat\g[[u,u^{-1}]]$). 

Recall that $r(z)=\frac{\Lambda}{z}+\text{regular part}$,
where $\Lambda$ is the inverse of the invariant form on $\g$.  
Therefore, $t(u)\Omega=\sum_a a\o a_+(u)\Omega$, where the summation 
is taken over an orthonormal basis of $\g$.  

Let $Y^0$ be the reduction of $Y$ mod $h$. We obtain 
$$
Y^0(t(u)\Omega,z)=t(u+z)-t^*(u+z)=
\Lambda\delta(x-u-z)=\sum_a a\o a(u+z).\tag 2.15
$$
In particular, $Y(a_+(u)\Omega,z)=a(u+z)$. 
This implies first part of the proposition.

Now let us compute the quasiclassical structure on $V^0$ 
which is induced by $V$. 
Let $\Cal S(z)=1+hs(z)+O(h^2)$. 
Then we obtain from the definition of $\Cal S(z)$:
$$
s^{34}(z)\circ t^{13}(u)t^{24}(v)=-[[r^{12}(u-v+z),t^{13}(u)]t^{24}(v)]
\tag 2.16
$$
By Proposition 1.8(ii), this implies that $s$ is given by (1.23), 
where $\rho=-r$.  
$\square$\enddemo

In view of Proposition 2.4, it is natural to call 
$V$ a quantum affine VOA. We will denote this quantum VOA by 
$V_q(sl_N,K,R)$. 

\head 3. The coinvariant construction 
of the quantum affine VOA for a rational R-matrix. \endhead

\subhead 3.1. The coinvariant construction of the affine VOA\endsubhead

There is a well known construction of the affine VOA
which does not require writing any formulas. This is the so-called 
coinvariant construction, which works as follows. 

As before, we let $V$ be the Weyl module over the affine 
Lie algebra $\hat\g$ with central charge $K$. Let $\Omega$
be a highest weight vector of $V$. 
 
Now we define 
the multiplication map $Y$ on $V$. We will first consider the case 
when $K$ is irrational. Consider 
the space $M=V\o V\o V^*$ (where $V^*$ is the graded dual to $V$). 
This space has an action of the Lie algebra 
$L=\hat\g\oplus\hat\g\oplus\hat\g/
((c,0,0)=(0,c,0)=(0,0,c))$ (here the action of $\hat\g$ on $V^*$ is as in the 
dual module to $V$,twisted by the transformation $t\to t^{-1}$). 
Let $L_{rat}\subset L$ be the Lie subalgebra consisting of triples 
of Laurent series which are expansions of a single rational 
$\g$-valued function $a(w)$ with no poles outside of $z,0,\infty$, 
with respect to the local parameters $w-z,w,1/w$. It can be shown 
that the space of invariant functionals $(M/LM)^*$ is 1-dimensional. Let 
$Y(z)\in (M/LM)^*$ be so normalized that $Y(z)(\Omega,*,*):V\o V^*\to k$ 
is the standard 
pairing. The element $Y(z)$ can be regarded as a map $V\o V\to \hat V$, 
where $\hat V$ is the completion of $V$ with respect to the grading. 
As a function of $z$, $Y(z)$ can be decomposed into a Laurent series, 
so $Y$ can be regarded as a map $V\o V\to V((z))$. 
This is the desired multiplication map. 

One can show that the map $Y$ extends by
continuity to all rational $K$. 
Further, it can be shown that there exists a unique $D:V\to V$
such that $V$ equipped with $(Y,D,\Omega)$ is a VOA, and that this VOA 
is isomorphic to $V(\g,K)$. Thus we have obtained an algebro-geometric 
construction of $V(\g,K)$. 

In the next sections we will generalize this construction to 
the quantum case for $\g=sl_N$. Since the construction uses 
a global object (projective line), our generalization works 
only for an R-matrix which is defined globally, i.e. the 
rational R-matrix. We do not know how to generalize this construction in 
the trigonometric and elliptic case. 

\subhead 3.2. Quantum current algebras with points at infinity
\endsubhead

{\bf 3.2.1.}
  From now on we let $R(u)$ be the normalized rational $R$-matrix 
$$
R(u)=f(u)(1-\frac{h(\sigma-1/N)}{N(u-h/N^2)}),\tag 3.1
$$
where $\sigma$ is the permutation, and $f(u)$ is the normalizing scalar 
function. 

Let us recall some constructions from \cite{EK4}. 
Let $F$ denote the algebra $F_0(R)$ defined 
in Chapter 2, and
$\hat F$ denote the corresponding central extension $\hat F_0(R)$. 
Let $U\supset F$ be the quantized 
universal enveloping algebra $U_0(R)$ corresponding to $F$ (see \cite{EK3}, 
Chapter 3). Let $\hat U$ be the central extension of $U$
(the quantized universal enveloping algebra corresponding to $\hat F$).

Let $F_\z$, $\hat F_\z$ be the corresponding 
factored Hopf algebras $F_0(R)_\z$, 
$\hat F_0(R)_\z$ ($\z=(z_1,...,z_N)\in \C^n[[h]]$).
Let $U_\z,\hat U_\z$ be the quantized universal 
enveloping algebras corresponding to $F_\z,\hat F_\z$. 

To proceed with our construction, we need to generalize 
the definition of $F_\z,\hat F_\z,U_\z,\hat U_\z$ as follows. 

Consider the Hopf algebra $F_*$ generated by a series
$T^*(v)=1+\sum_{j\ge 1}T^*_{-j}v^{-j}$, satisfying the defining 
relations 
$$
R^{12}(u-v)T^{*13}(u)T^{*23}(v)=T^{*23}(v)T^{*13}(u)R^{12}(u-v),
\qdet(T^*(u))=1,\tag 3.2
$$
with coproduct, counit, and antipode given by 
the usual formulas
$$
\Delta(T^*(u))=T^{*12}(u)T^{*13}(u), \e(T^*(u))=1, S(T^*(u))=T^*(u)^{-1}.\tag 
3.3
$$
Let $U_*\supset F_*$ be the quantized universal enveloping algebra 
corresponding to $F_*$. 

The Hopf algebras $F,F_*$ are deformations of the function algebras 
on the groups $SL_n[[t]]$, $1+t^{-1}sl_N[[t^{-1}]]$, respectively. 
The Hopf algebras algebras $U$, $U_*$ are isomorphic to the dual Yangian 
(with opposite coproduct) and the ordinary 
Yangian for $sl_N$, respectively.  

{\bf Remark.} The flatness of the algebra defined by (3.2) 
depends on the fact that $R$ is the rational R-matrix. 

There is a natural pairing $(,)$ between $F$ and $F_*$, 
induced by the 
Hopf algebra embedding 
$F_*\to F^{*op}$, which is defined in section 4.1 of \cite{EK4}. 
It is defined by the formula $(T^{*13}(u),T^{23}(v))=R^{12}(u-v)$. 

Define the algebra $F_{\z,\infty}$ generated by 
$F_\z$ and $F_*$ with commutation relations
$$
R^{12}(u-v+z_i)T^{13}_i(u)T^{*23}(v)=T^{*23}(v)T^{13}_i(u)R^{12}(u-v+z_i),
\tag 3.4
$$
It is clear that $F_{\z,\infty}$ is a factored 
algebra in the sense of \cite{EK3}, with factors 
$F_\z$ and $F_*$. It is also a factored Hopf algebra 
with coproduct, counit and antipode induced from the factors. 
One can also define the corresponding quantum universal enveloping algebra 
$U_{\z,\infty}$. 

Similarly one defines the central extension $\hat F_{\z,\infty}$.
Let $\hat F_*=F_*\o k[c]$ as an algebra, with coproduct 
$$
\Delta(T^*(u))=T^{*12}(u-hc^3/2)T^{*13}(u+hc^2/2), 
\Delta(c)=c\o 1+1\o c.\tag 3.5
$$
Define the algebra $\hat F_{\z,\infty}$ to be the algebra generated by 
$\hat F_\z$ and $\hat F_*$ with commutation relations
$$
 R^{12}(u-v+z_i-\frac{h}{2}(c_i-c))
T^{13}_i(u)T^{*23}(v)=T^{*23}(v)T^{13}_i(u)
 R^{12}(u-v+z_i+\frac{h}{2}(c_i-c)),
\tag 3.6
$$
and the condition that $c_i$ and $c$ are central in the whole algebra. 
It is clear that $\hat F_{\z,\infty}$ is a factored 
algebra in the sense of \cite{EK3}, with factors 
$\hat F_\z$ and $\hat F_*$. It is also a factored Hopf algebra 
with coproduct, counit and antipode induced from the factors. 
Let $\hat U_{\z,\infty}$ denote the corresponding quantum universal 
enveloping algebra. 

The algebra $U_*$, which is, as we mentioned, the Yangian of $sl_N$, 
contains a copy of $U(sl_N)[[h]]$ as a Hopf subalgebra. 
This Hopf subalgebra is generated by the components of 
$T_{-1}^*$, which span the Lie algebra $sl_N$ of primitive elements. 

{\bf 3.2.2.}
In \cite{EK4}, section 3.4, 
for any finite-dimensional $F$-comodule $W$,
we constructed local dimodules $W_K(z_i)$. 
They are $(\hat F_\z,\hat F_\z,C)$-dimodules, where $C$ is the ideal 
generated by $c_1,...,c_n$. It is easy to see that they are also 
$(\hat U_\z,\hat U_\z,C)$-dimodules. 

It is necessary for us to extend 
$\hat W_K(z_i)$ to a $(\hat U_{\z,\infty},\hat U_{\z,\infty},\tilde C)$-
dimodule, where $\tilde C$ is the ideal generated by $c_1,...,c_n,c$. 
To do this, it is enough to define an action of $T^*(u)$ and $c$ 
on $\hat W_K(z_i)$. We do this by analogy with formulas (3.17), (3.21) 
of \cite{EK4}: 
$$
T^*(u)\to (T^{*13}(u-z_i),T^{24}_{\hat W_K})_{34},\ c=K,\tag 3.7
$$
where $(,):F_*\times F\to k[[h]]$ is the natural pairing.

Also, we need to define local dimodules
concentrated at infinity. For this purpose define 
the subalgebras $\bar U_\z$ 
of $\hat U_{\z,\infty}$ generated by $\hat U_\z$, the matrix elements of
$T_{-1}^*$, and the element $c$. As a vector space, this algebra 
is naturally isomorphic to $\hat U_{\z}\o U(sl_N)\o k[c]$, because of the 
commutation relation
$$
[T_{-1}^{*23},T^{13}_i(u)]=-\frac{h}{N}[\sigma^{12},T^{13}_i(u)].
$$ 
 
Let $\bold 1_K(\infty)$ 
be the 1-dimensional
$(\bar U_\z,\hat U_{\z,\infty},C)$-dimodule, which is trivial as a 
comodule and as an $\bar U_\z$-module, and such that $c_i=c=K$. 
Define
$$
\hat {\bold 1}_K(\infty)=
\text{Ind}_{(\bar U_\z,\hat U_{\z,\infty},C)}
^{(\hat U_{\z,\infty},\hat U_{\z,\infty},C)}\bold 1_K(\infty).\tag 3.8
$$
  
Let $Z=U_*\o_{U(sl_N)}k$ be an induced Yangian module. 
Observe that as a vector space, $\hat W_K(z_i)$ is naturally 
identified with $W\o U$, and 
$\hat {\bold 1}_K(\infty)$ is naturally identified with $Z$. 

{\bf 3.2.3.}
Now by analogy with Section 3.5 of 
\cite{EK4} define a  global dimodule $M_K(W^1,...,W^n,\z,\infty)=
\hat W_K^1(z_1)\o...\o \hat W_K^n(z_n)\o \hat{\bold 1}_K(\infty)$. 
It is a $(\hat U_{\z,\infty},\hat U_{\z,\infty},C)$-dimodule. 

Let $\tilde U_{\z,\infty}$ be the quotient of $\hat U_{\z,\infty}$ 
by the ideal generated by the relations $c_i=c,i=1,...,n$. It is clear that 
$\tilde U_{\z,\infty}$ is a Hopf algebra. 
An important property of $\tilde U_{\z,\infty}$, which 
$\hat U_{\z,\infty}$ does not have, is that $U_{\z,\infty}$ is a subalgebra
of $\tilde U_{\z,\infty}$.

It is clear that the dimodule $M_K(W^1,...,W^n,\z,\infty)$
descends to a $(\tilde U_{\z,\infty},\tilde U_{\z,\infty},\<c\>)$-dimodule. 
For brevity we will denote this dimodule $M_K(\z,\infty)$. 

{\bf 3.2.4.}
Since $U_{\z,\infty}\subset \tilde U_{\z,\infty}$ is a subalgebra, 
we have an action of $U_{\z,\infty}$ on 
$M_K(\z)$. Define the space of invariant
functionals 
$$
B_K(W^1,\dots ,W^n,\z,\infty)=B_K(\z,\infty):=
\Hom_{U_{\z,\infty}}(M_K(\z,\infty),k[[h]]).\tag 3.9
$$

The elements of $B_K(\z,\infty)$ are quantum analogues 
of conformal blocks in the Wess-Zumino-Witten model.
Therefore we will call them ``quantum conformal blocks''.

We have a natural evaluation map $\xi:B_K(\z,\infty)\to (W^1\o...\o W^n)^*$ 
defined by
$$
\xi(f)(\v_1\o...\o \v_n)=f(\v_1\o...\o \v_n\o 1), \v_i\in W^i(z_i)\subset 
\hat W^i_K(z_i).
$$

Similarly to Proposition 3.5 of \cite{EK4}, 
 the map $\xi$ is a linear isomorphism.

\subhead 3.3. The quantum coinvariant construction\endsubhead

{\bf 3.3.1.}
We will be most interested in the case when $W_i$ are all equal to the 
trivial module. In this case the space $M_K(\z,\infty)$ 
is isomorphic to $U^{\o n}\o Z$, and a conformal block is unique 
up to scaling. 

For $n=0$, this conformal block $b_1$ 
is just the descent of the counit of $U_*$
to $Z$. 

For $n=1$, the conformal block $b_2$ is a bilinear form 
$U\o Z\to k[[h]]$. In the classical limit, $U$ and $Z$ can be 
naturally identified 
with the Weyl module $V_{0,K}$ of highest weight zero 
and central charge $K$ 
over the affine Kac-Moody algebra $\widehat{sl_N}$. 
In this case, $b_2$ corresponds to the contravariant 
(Shapovalov) form. In particular, if $K$ is irrational, 
this form is nondegenerate. 
 
Now consider the next special case 
$n=2$. By shifting of the coordinate it reduces to $\z=(z,0)$, where $z$ 
is a nonzero point. In this case, $M_K(z,0,\infty)=U\o U\o Z$, 
and the unique up to scaling conformal block will be denoted by $b_3(z)$
(we normalize $b_3(z)$ so that $b_3(z)(1,1,1)=1$). 
We can regard $b_3(z)$ as a map $U\o U\o Z\to k[z,z^{-1}][[h]]$. 

We assume that $K$ is irrational. In this case 
there exists a unique map $Y(z):U\o U\to U((z))$ defined by the equation
$b_2(Y(z)(\u_1\o \u_2),\v)=b_3(z)(\u_1,\u_2,\v)$. Indeed, 
this statement is well known classically, and the quantum statement 
is proved by perturbation argument. 

{\bf 3.3.2.}
Now let us construct a quantum deformation of the affine VOA. 
We have to produce a set $(V,Y,T,\Omega,\Cal S)$ satisfying
some conditions. The first two elements have already been constructed: 
we set $V=U$ and $Y:V\o V\to V((z))$ as above.  
Thus, to finish the construction, we need to 
define elements $\Cal S,T,\Omega$ and prove their properties. 

We start with $\Cal S$. By definition, $\Cal S(w)$ should be a map 
$V\o V\to V\o V\o k((w))$.

Consider the 
Drinfeld pseudotriangular structure $\RR(u)\in F^*\o F^*((u))$. 
Set $\Cal S=\RR^{-1}$ (this makes sense since $V=\hat{\bold 1}_K$ has an 
dimodule structure, in particular the structure of a module over $F^*$). 

Further, we set $\Omega=1$. 
Finally, we 
define the Sugawara operator $T$ to be the operator $D$
given by $e^{zD}(T_1(u_1)...T_n(u_n)\Omega)=T_1(u_1+z)...T_n(u_n+z)\Omega$. 

\proclaim{Theorem 3.1} $(V,Y,D,\Omega,\Cal S)$ is a 
quantum vertex operator algebra. 
\endproclaim

\demo{Proof} Properties of $\Cal S$ are easy to prove. 

Axiom (QA1) follows from the fact that the operator product 
$Y(\u_1,z_1)Y(\u_2,z_2)\w$ satisfies the equation
$b_2(Y(\u_1,z_1)Y(\u_2,z_2)\w,\v)=b_4(z_1,z_2)(\u_1,\u_2,\w,\v)$,
where $b_4$ is the unique
element of $B_K(z_1,z_2,0,\infty)$ such that $b_4(z_1,z_2)(1,1,1,1)=1$, 
and the fact that 
$\sigma\RR(z_1-z_2)$ is an isomorphism of dimodules 
$\hat V^1_K(z_1)\o\hat V^2_K(z_2)\to \hat V^2_K(z_2)\o\hat V^1_K(z_1)$. 

Axioms (QA2) and (QA3) are deduced
from the definition of $Y$, $D$ 
and $\Omega$ as in the classical case.

It remains to prove axiom (QA4). 
As we know from Proposition 1.11(ii), 
if $V/hV$ is nondegenerate, then (QA4) is equivalent to 
associativity. But $V/hV$ is the affine VOA, which is
nondegenerate for irrational $K$ by Proposition 1.12. 
Thus, it suffices to prove the associativity 
identity (1.12). This identity follows from the fact that 
both sides of (1.12) are elements of $B_K(z,w,0,\infty)$ satisfying the same 
normalization condition. 
$\square$\enddemo
 
We will denote the obtained quantum VOA by $\tilde V_q(sl_N,K,R)$. 

\subhead 3.4. Explicit computation of $Y$\endsubhead

{\bf 3.4.1.}
Here we compute the structure of $\tilde V_q(sl_N,K,R)$ completely
in terms of the generators $T(u),T^*(u)$ of the double Yangian, and show 
that it coincides with $V_q(sl_N,K,R)$. 

\proclaim{Proposition 3.2} The map $Y$
for $\tilde V(sl_N,K,R)$ is defined 
by the formula 
$$
\gather
Y(T^{1,n+1}(u_1)...T^{n,n+1}(u_n)\Omega,z)=\\
T^{1,n+1}(u_1+z)...T^{n,n+1}(u_n+z)T^{*n,n+1}(u_n+z+Kh/2)...
T^{*1,n+1}(u_1+z+Kh/2).
\tag 3.10
\endgather
$$
\endproclaim

\demo{Proof} 
It is clear that both sides of (3.10) coincide on $\Omega$, so it is enough 
to show that they commute in the same way with $T(v)$. 

We first note that $T$ and $T^*$ satisfy the relation
$$
 R^{12}(u-v-Kh/2)T^{13}(u)T^{*23}(v)=T^{*23}(v)T^{13}(u)
 R^{12}(u-v+Kh/2).
\tag 3.11
$$
This relation follows from the definition of the dimodule $\hat\bold 1_K$. 

Denote the right hand side of (3.10) by $\bold T_n(u_1,...,u_n,z)$. 
Then a direct computation using (3.11) yields 
$$
\gather
T^{0,n+1}(v)\bold T_n^{1...n,n+1}
(u_1,...,u_n,z)=\\  R^{01}(v-z-u_1)^{-1}... R^{0n}(v-z-u_n)^{-1}
\bold T_n^{1...n,n+1}(u_1,...,u_n,z)\times\\ 
 R^{0n}(v-z-u_n-Kh)... R^{01}(v-z-u_1-Kh)
T^{0,n+1}(v).\tag 3.12\endgather
$$

On the other hand, 
since $Y$ was constructed 
from a coinvariant, we have the identity
$$
T^{01}(v)Y(z)=Y(z)T^{*01}(v-z-Kh/2)T^{02}(v).\tag 3.13
$$
Using this relation and (3.11), we get that
(3.12) remains valid if $\bold T_n(u_1,...,u_n,z)$ is replaced with 
 $Y(T^{1,n+1}(u_1)...T^{n,n+1}(u_n)\Omega,z)$.
The proposition is proved. 
$\square$\enddemo

\proclaim{Corollary 3.3}
The quantum VOA $\tilde V_q(sl_N,K,R)$ coincides with  
$V_q(sl_N,K,R)$.
\endproclaim

\demo{Proof} First of all, formula (3.11) implies that the operator 
series $T^*(u)$ defined in this chapter is the same as in chapter 2. 
Therefore, the multiplication map in both algebras is the same
by Proposition 3.2. 

Since $K$ is irrational, the classical limit of 
$V_q$ and $\tilde V_q$ is a nondegenerate VOA. 
But if $V_1,V_2$ are any quantum VOA realized on the same space $V$, such that 
$Y_{V_1}=Y_{V_2}$, and the classical limits of $V_1,V_2$ are nondegenerate, 
then $\Cal S_{V_1}=\Cal S_{V_2}$ and hence $V_1=V_2$. 
This implies the corollary.   
$\square$\enddemo

\proclaim{Corollary 3.4} The construction of section 3.3 extends by continuity 
to rational values of $K$. 
\endproclaim

\demo{Proof} Clear from Corollary 3.3.$\square$\enddemo

{\bf 3.4.2.} 
Recall that 
the Sugawara operator of the affine VOA 
$V(\g,K)$ at the non-critical value of 
the central charge $K$ (i.e. $K\ne -g$, where $g$ is the dual 
Coxeter number of $\g$) is given by the following explicit formula, 
called the Sugawara construction:
$$
D=-\frac{1}{2(K+g)}\text{Res}_{z=0}\sum_a:a(z)^2:dz.\tag 3.14
$$
This construction can be generalized to the quantum case. 

\proclaim{Proposition 3.5} If $K\ne -N$ then 
the Sugawara element $D$ of the quantum VOA 
$\tilde V_q(sl_N,K,R)$ is given by the formula 
$$
D=-\frac{1}{K+N}\ln Q,\tag 3.15
$$
 where $Q$ is the 
quantum Sugawara element defined in \cite{EK4}. 
\endproclaim

\demo{Proof} This follows from formula (5.3) in \cite{EK4}
$\square$\enddemo

\end